\documentclass[10pt]{amsart}

\usepackage{graphicx,amscd,amssymb,verbatim}
\usepackage[dvips]{hyperref}
\usepackage[leqno]{amsmath}

 \newtheorem{thm}{Theorem}[section]
 \newtheorem{Corollary}{Corollary}[section]

 \numberwithin{equation}{section}
 \def\Gr{\operatorname{Gr}}

 \DeclareMathOperator{\Id}{Id}
 
 \newcommand{\To}{\longrightarrow}

 \newcommand{\norm}[1]{\left\Vert#1\right\Vert}

\begin{document}

\title{Four-dimensional Riemannian manifolds with commuting higher order Jacobi operators}
\author{Zhivko Zhelev, Maria Ivanova${}^1$ and Veselin Videv}

\address{Zh. Zhelev: Department of Mathematics and Informatics,
 Trakia University, St. Zagora,
 6000 Bulgaria. Email: \tt{zhelev@uni-sz.bg}}

\address{V. Videv: Department of Mathematics and Informatics,
Trakia University, St. Zagora,
 6000 Bulgaria. Email: \tt{videv@uni-sz.bg}}

\address{M. Ivanova: Department of Mathematics and Informatics,
Sofia University, Sofia, 1000 Bulgaria. Email:
\tt{m\_stoeva@abv.bg}}

\date{}

\begin{abstract}
We consider four-dimensional Riemannian manifolds with commuting
higher order Jacobi operators defined on two-dimensional
orthogonal subspaces (polygons) and on their orthogonal subspaces.\\
\indent More precisely, we discuss higher order Jacobi operator
$\mathcal{J}(X)$ and its commuting associated operator
$\mathcal{J}(X^{\perp})$ induced by the orthogonal complement
$X^{\perp}$ of the vector $X$, i. e.
$\mathcal{J}(X)\circ\mathcal{J}(X^{\perp})=\mathcal{J}(X^{\perp})\circ \mathcal{J}(X)$.\\
\indent At the end some new central theorems have been cited. The
latter are due to P. Gilkey, E. Puffini and V. Videv, and have
been recently obtained.
\end{abstract}

\keywords{Einstein manifold, Higher order Jacobi operator, Jacobi
operator, Ricci operator.
${}^1$Corresponding author}

 \maketitle

\section{Preliminaries} Let $(M,g)$ be a $n$-dimensional
Riemannian manifold with a metric tensor $g$. Tangent space at a
point $p\in M$ we denote by $M_p$, and let $S_pM$ be the set of
unit vectors in $M_p$,\\
 i. e. $S_p(M):=\{z\in
M_p\,|\,\,\norm{g(z,z)}=1\}$. Let $\mathcal{F}(M)$ be the algebra
of all smooth functions on $M$ and $\mathcal{X}(M)$ be the
$\mathcal{F}(M)$-module of all smooth vector fields over $M$. Let
also
$$
R(X,Y)Z=\nabla_X\nabla_YZ-\nabla_Y\nabla_XZ-\nabla_{[X,Y]}Z,\qquad
X,Y,Z\in\mathcal{X}(M)
$$
be the $(1,3)$ curvature tensor of the Levi-Civita connection
$\nabla$. We define
$$
R(X,Y,Z,U):=g(R(X,Y,Z),U)
$$
to be the associated $(0,4)$-curvature tensor which satisfied the
following algebraic properties:
$$
\begin{array}{l}
   \mbox{i)}\quad R(X,Y,Z,U)=-R(Y,X,Z,U),\\
    \mbox{ii)}\quad R(X,Y,Z,U)=-R(X,Y,U,Z),\\
   \mbox{iii)}\quad R(X,Y,Z,U)+R(Y,Z,X,U)+R(Z,X,Y,U)=0\,\,\mbox{(first Bianchi
   identity)},\\
   \mbox{iv)}\quad R(X,Y,Z,U)=R(Z,U,X,Y).
 \end{array}
$$
\indent In the Riemannian geometry the following differential
equality is also true:
$$
\begin{array}{l}
\mbox{v)}\quad
(\nabla_XR)(Y,Z,W)+(\nabla_YR)(Z,X,W)+(\nabla_ZR)(X,Y,W)=0\,\,\mbox{(second
Bianchi identity)},
\end{array}
$$
where
$$
(\nabla_XR)(Y,Z,W):=\nabla_X(R(Y,Z)W)-R(\nabla_XY,Z)W-R(Y,\nabla_XZ)W-R(Y,Z)\nabla_XW,
$$
and $\nabla_XR$ is the covariant derivative of the
$(1,3)$-curvature tensor $R$ with respect to $X$,
$X,Y,Z\in\mathcal{X}(M)$.\\
 \indent Let $\mathcal{J}\colon M_p\To M_p$ be the Jacobi
operator defined by:
\begin{equation}
\mathcal{J}(X)U=R(U,X,X).
\end{equation}
\indent One can easily see that $\mathcal{J}(X)X=0$ and
$g(\mathcal{J}(X)Y,Z)=g(Y,\mathcal{J}(X)Z)$ which means that
Jacobi operator is a symmetric linear operator.\\
\indent Jacobi operator can be diagonalized in the Riemannian
geometry. In this case we say that $g$ is \textit{Osserman metric}
if the eigenvalues of the Jacobi operator are constant over the
tangent bundle $S(M):=\bigcup_{p\in M}S_pM$. If $(M,g)$ is a rank
one locally symmetric space, i. e. $\nabla R=0$, where $\nabla$ is
the connection with all positive or all negative sectional
curvatures \cite{Wang} or $(M,g)$ is flat, i. e. $R=0$, the group
of local isometries acts transitively on $S(M)$ and each Jacobi
operator has constant eigenvalues. Osserman \cite{Oss} conjectured
that the opposite is also true and this was confirmed by Chi when
$n=4$ and $n\equiv (2\mod 4)$ \cite{Chi} and by Nikolaevsky
\cite{Nik} when $n\neq
16$.\\
\indent Gilkey, Stanilov and Videv \cite{Vid1} introduced a new
operator which they called \textit{general Jacobi operator of
order $k$ or $k$-order Jacobi operator}. More precisely, if
$\{Y_i\}_{i=1}^k$ is any orthonormal basis for an arbitrary
$k$-plane $\pi\in M_p$, the higher $k$-order Jacobi operator is
defined by:
\begin{equation}
\mathcal{J}(\pi)\colon Y\To\sum_{1\le i\le
k}R(Y,Y_i)Y_i=\sum_{1\le i\le k}\mathcal{J}(Y_i)Y.
\end{equation}
\indent It can be easily verified that this operator does not
depend on the basis of $\pi$.

\section{Some commutativity conditions}
\indent Another variety of problems, connected to the higher order
Jacobi operator, emerged thanks to Stanilov and Videv \cite{Vid2}.
They are connected with some commutativity conditions forced on
(1.2). Recently Brozos-V\'{a}zquez and Gilkey \cite{Bro} were able
to prove the following
\begin{thm}
Let $(M,g)$ be a Riemannian manifold, $\dim M\ge 3$. Then
\begin{itemize}
\item[(A)] $(M,g)$ is \textit{flat} iff
$\mathcal{J}(X)\mathcal{J}(Y)=\mathcal{J}(Y)\mathcal{J}(X)$ for
arbitrary vectors $X,Y\in M_p$; \item[(B)] $(M,g)$ is a
\textit{manifold with a constant sectional curvature} iff
$\mathcal{J}(X)\mathcal{J}(Y)=\mathcal{J}(Y)\mathcal{J}(X)$ for
arbitrary vectors $X,Y\in M_p$ such that $X\perp Y$.
\end{itemize}
\end{thm}
\indent In this paper authors will characterize indecomposible four-dimensional
Riemannian manifolds that satisfy the following two
conditions:\\

 \textit{ For arbitrary unit vector $X\in M_p$, $p\in
M$, we have:
\begin{itemize}
\item[{\bf{(C1)}}]$\mathcal{J}(X)\circ\mathcal{J}(X^{\perp})=\mathcal{J}(X^{\perp})\circ
\mathcal{J}(X)$, where $X^{\perp}$ is the orthogonal complement of
$X$ in $M_p$.
\end{itemize}
For arbitrary $2$-plane $\alpha\subset M_p$, $p\in M$, we have
\begin{itemize}
\item[{\bf{(C2)}}] $\mathcal{J}(\alpha
)\circ\mathcal{J}(\alpha^{\perp})=
\mathcal{J}(\alpha^{\perp})\circ\mathcal{J}(\alpha )$, where
$\alpha^{\perp}$ is the orthogonal complement of $\alpha$ in
$M_p$.
\end{itemize}
}
 \indent Our main goal is to prove the following
 \begin{thm}
 Let $(M,g)$ be a four-dimensional indecomposible Riemannian manifold. Then the
 following are equivalent:
 \begin{itemize}
 \item[(a)]  Equality {\bf{(C1)}} holds for arbitrary unit vector
 $X\in M_p$, $p\in M$;
 \item[(b)] Equality {\bf{(C2)}} holds for arbitrary $2$-plane
 $\alpha\subset M_p$, $p\in M_p$;
 \item[(c)] $(M,g)$ is Einstein.
 \end{itemize}
 \end{thm}

 \textit{Proof.} $\bf{(a)\Longrightarrow (c)}$ Let
 $\{e_1,e_2,e_3,e_4\}$ be an orthonormal basis for $M_p$, $p\in
 M$. Curvature operator matrices
 $\mathcal{J}_{\{e_1,e_2,e_3\}}$ and $\mathcal{J}_{\{e_4\}}$ then have
 the form:
 \begin{equation}
  \begin{pmatrix}
  K_{12}+K_{13} & R_{1332} & R_{1223} & \rho_{14}\\
  R_{1332} & K_{12}+K_{13} & R_{2113} & \rho_{24}\\
  R_{1223} & R_{2113} & K_{13}+K_{23} & \rho_{34}\\
  \rho_{14} & \rho_{24} & \rho_{34} & K_{14}+K_{24}+K_{34}
  \end{pmatrix},
 \end{equation}
and
\begin{equation}
  \begin{pmatrix}
    K_{14} & R_{1442} & R_{1443} & 0\\
    R_{1442} & K_{24} & R_{2443} & 0\\
    R_{1443} & R_{2443}& K_{34} & 0\\
    0 & 0 & 0 & 0
    \end{pmatrix},
  \end{equation}
  where $\rho_{ij}=\rho
  (e_i,e_j):=\displaystyle{\sum_{k=1}^{4}g(R(e_k,e_i)e_j,e_k)}$
  are the components of the Ricci $(0,2)$-tensor $\rho$
  and $K_{ij}:=g(R(e_i,e_j)e_j,e_i)$, $i,j=1,\ldots ,4$.\\
  \indent We have the matrix equality
  \begin{equation}
   \mathcal{J}_{\{e_1,e_2,e_3\}}\circ\mathcal{J}_{\{e_4\}}=\mathcal{J}_{\{e_4\}}\circ\mathcal{J}_{\{e_1,e_2,e_3\}},
   \end{equation}
  which leads us to the equations:
  \medskip
  \begin{itemize}
  \item[($e_1)$]$\quad$
  $K_{14}(R_{1224}+R_{1334})+R_{1442}(R_{2114}+R_{2334})+R_{1443}(R_{3114}+R_{3224})=0$,
  \item[$(e_2)$]$\quad$
 $ R_{1442}(R_{1224}+R_{1334})+K_{24}(R_{2114}+R_{2334})+R_{2443}(R_{3114}+R_{3224})=0$,
  \item[$(e_3)$]$\quad$
  $R_{1443}(R_{1224}+R_{1334})+R_{2443}(R_{2114}+R_{2334})+K_{34}(R_{3114}+R_{3224})=0$.
  \end{itemize}
  \medskip

  \indent We do a cyclic change of indecies $1\To 2\To 3\To 4\To
  1$ in $(e_1)$, $(e_2)$, $(e_3)$, and get the equations
  \medskip

\begin{itemize}
  \item[($e_1^1)$]$\quad$
  $K_{12}(R_{1332}+R_{1442})+R_{2113}(R_{1223}+R_{1443})+R_{2114}(R_{1224}+R_{1334})=0$,
  \item[$(e_1^2)$]$\quad$
 $ R_{2113}(R_{1332}+R_{1442})+K_{13}(R_{1223}+R_{1443})+R_{3114}(R_{1224}+R_{1334})=0$,
  \item[$(e_1^3)$]$\quad$
  $R_{2114}(R_{1332}+R_{1442})+R_{3114}(R_{1223}+R_{1443})+K_{14}(R_{1224}+R_{1334})=0$.
  \end{itemize}
\medskip

  \indent Doing a cyclic change of indecies $1\To 2\To 3\To 4\To
  1$ in $(e_1^1)$, $(e^2_1)$ and $(e_3^1)$, we get
\medskip

\begin{itemize}
  \item[($e_2^1)$]$\quad$
  $K_{23}(R_{2113}+R_{2443})+R_{3224}(R_{2114}+R_{2334})+R_{1223}(R_{1332}+R_{1442})=0$,
  \item[$(e_2^2)$]$\quad$
 $ R_{3224}(R_{2113}+R_{2443})+K_{24}(R_{2114}+R_{2334})+R_{1224}(R_{1332}+R_{1442})=0$,
  \item[$(e_2^3)$]$\quad$
  $R_{1223}(R_{2113}+R_{2443})+R_{1224}(R_{2114}+R_{2334})+K_{12}(R_{1332}+R_{1442})=0$.
  \end{itemize}
  \medskip

\indent Another cycling change $1\To 2\To 3\To 4\To
  1$ in $(e_2^1)$, $(e_2^2)$ and $(e_2^3)$ will give us
\medskip

\begin{itemize}
  \item[($e_3^1)$]$\quad$
  $K_{34}(R_{3114}+R_{3224})+R_{1334}(R_{1223}+R_{1443})+R_{2334}(R_{2113}+R_{2443})=0$,
  \item[$(e_3^2)$]$\quad$
 $ R_{1334}(R_{3114}+R_{3224})+K_{13}(R_{1223}+R_{1443})+R_{1332}(R_{2113}+R_{2443})=0$,
  \item[$(e_3^3)$]$\quad$
  $R_{2334}(R_{3114}+R_{3224})+R_{1332}(R_{1223}+R_{1443})+K_{23}(R_{2113}+R_{2443})=0$.
  \end{itemize}
\medskip

  \indent Solving $(e_1)$, $(e_2)$, $(e_3)$ with respect to
  $R_{1224}+R_{1334}$, $R_{2114}+ R_{2334}$, $R_{3114}+R_{3224}$,
  using \textit{Maple},
   we get:

  $$
  R_{1224}+R_{1334}=R_{2114}+R_{2334}=R_{3114}+R_{3224}=0,
  $$
\medskip

\noindent since the above is in fact the trivial solution to the
system of equations $(e_1)$, $(e_2)$, $(e_3)$ which is
homogeneous:
$$
\begin{pmatrix}
K_{12} & R_{1442} & R_{1443}\\
R_{1442} & K_{24} & R_{2443}\\
R_{1443} & R_{2443} & K_{34}
\end{pmatrix}
\begin{pmatrix}
R_{1224}+R_{1334}\\
R_{2114}+R_{2334}\\
R_{3114}+R_{3224}
\end{pmatrix}
=
\begin{pmatrix}
0\\
0\\
0
\end{pmatrix}.
$$
\indent Analogously, solving the other two homogeneous systems
$(e_1^1)$, $(e_1^2)$, $(e_1^3)$ and $(e_3^1)$, $(e_3^2)$,
$(e_3^3)$ we get that

$$
R_{1332}+R_{1442}=R_{1223}+R_{1443}=R_{1224}+R_{1334}=0,
$$

and
$$
R_{1223}+R_{1443}=R_{2113}+R_{2443}=R_{3114}+R_{3224}=0.
$$
\indent From (2.3) we also get
\medskip

 \begin{itemize}
  \item[($e_4)$]$\quad$
  $(K_{14}-K_{24})R_{1332}+R_{2113}R_{1443}-R_{2443}R_{1223}=0$,
  \item[$(e_5)$]$\quad$
 $ (K_{14}-K_{34})R_{1223}+R_{2113}R_{1442}-R_{2443}R_{1332}=0$,
  \item[$(e_6)$]$\quad$
  $(K_{24}-K_{34})R_{2113}+(K_{13}-K_{12})R_{2443}+R_{1223}R_{1442}-R_{1332}R_{1443}=0$.
  \end{itemize}
\medskip

  \indent By a cycling change of indecies $1\To 2\To 3\To 4\To
  1$ in $(e_4)$, $(e_5)$ and $(e_6)$, we get
\medskip

  \begin{itemize}
  \item[($e_4^1)$]$\quad$
  $(K_{12}-K_{13})R_{2443}+R_{2114}R_{3224}-R_{3114}R_{2334}=0$,
  \item[$(e_4^2)$]$\quad$
 $ (K_{12}-K_{14})R_{2334}+R_{3224}R_{2113}-R_{3114}R_{2443}=0$,
  \item[$(e_4^3)$]$\quad$
  $(K_{13}-K_{14})R_{3224}+(K_{24}-K_{23})R_{3114}+R_{2334}R_{2113}-R_{2443}R_{2114}=0$.
  \end{itemize}
\medskip

\indent Repeating the same procedure two more times, we get
\medskip

\begin{itemize}

  \item[($e_5^1)$]$\quad$
  $(K_{23}-K_{24})R_{3114}+R_{1334}R_{1223}-R_{1224}R_{1443}=0$,
  \item[$(e_5^2)$]$\quad$
 $ (K_{23}-K_{12})R_{1443}+R_{3224}R_{1334}-R_{1224}R_{3114}=0$,
  \item[$(e_5^3)$]$\quad$
  $(K_{24}-K_{12})R_{1334}+(K_{13}-K_{34})R_{1224}+R_{1443}R_{3224}-R_{3114}R_{1223}=0$.
  \end{itemize}

  and

\begin{itemize}
 \item[($e_6^1)$]$\quad$
  $(K_{34}-K_{13})R_{1224}+R_{2334}R_{1442}-R_{1332}R_{2114}=0$,
  \item[$(e_6^2)$]$\quad$
 $ (K_{34}-K_{23})R_{2114}+R_{1334}R_{1442}-R_{1332}R_{1224}=0$,
  \item[$(e_6^3)$]$\quad$
  $(K_{13}-K_{23})R_{1442}+(K_{24}-K_{14})R_{1332}+R_{1334}R_{2114}-R_{2334}R_{1224}=0$.
  \end{itemize}

\medskip
 Further, from $(e_4)$, $(e_5)$, $(e_6)$; $(e_4^1)$,
$(e_4^2)$, $(e_4^3)$; $(e_5^1)$, $(e_5^2)$, $(e_5^3)$; $(e_6^1)$,
$(e_6^2)$, $(e_6^3)$ and using that

$$
\begin{array}{cc}
R_{2113}+R_{2443}=0, & R_{1332}+R_{1442}=0\\
R_{1223}+R_{1443}=0, & R_{1224}+R_{1334}=0\\
R_{2114}+R_{2334}=0, & R_{3114}+R_{3224}=0,
\end{array}
$$

we get the system of equations
$$\left|
\begin{array}{c}
K_{12}=K_{34}\\
\\
 K_{13}=K_{24}\\
 \\
 K_{14}=K_{23}\end{array}, \right.
 $$
 with respect to the basis $\{e_1, e_2, e_3, e_4\}$. The latter is
 equivalent to $(M,g)$ being an Einstein\cite{Sin}.$\hfill\square$\\

 $\bf{(c)\Longrightarrow (a)}$ Suppose $(M,g)$ is a four-dimensional
 Einstein manifold and let $X\in M_p$, $p\in M$ and $X^{\perp}$ is
 the orthogonal complement of $X$. Then $\rho
 =\lambda\Id,\,\,\lambda=\mbox{const.}$, and hence
 $$
 \begin{array}{c}
 \mathcal{J}(X)\circ\mathcal{J}(X^{\perp})-\mathcal{J}(X^{\perp})\circ\mathcal{J}(X)=
\mathcal{J}(X)\circ\mathcal{J}(X^{\perp})+\mathcal{J}(X)\circ\mathcal{J}(X)-\mathcal{J}(X)\circ\mathcal{J}(X)
-\mathcal{J}(X^{\perp})\circ\mathcal{J}(X)=\\
\\
\mathcal{J}(X)\circ\underbrace{\left
[\mathcal{J}(X^{\perp})+\mathcal{J}(X)\right ]}_{\rho}-
\underbrace{\left [\mathcal{J}(X)+\mathcal{J}(X^{\perp})\right
]}_{\rho}\circ\mathcal{J}(X)=\mathcal{J}(X)\circ\rho
-\rho\circ\mathcal{J}(X)=\\
\lambda\left (\mathcal{J}(X)\circ\Id-\Id\circ\mathcal{J}(X)\right
)=0.\\
\hfill\square
\end{array}
$$

\indent Analogously, one can prove, using \cite{Schouten}, the
following
\begin{Corollary}
Let $(M,g)$ be a three-dimensional Riemannian manifold. Then the
next two conditions are equivalent:
\begin{itemize}
\item[(i)]
$\mathcal{J}(X)\circ\mathcal{J}(X^{\perp})=\mathcal{J}(X^{\perp})\circ\mathcal{J}(X)$
for arbitrary $X,\,X^{\perp}\in M_p$, $p\in M$. \item[(ii)]
$(M,g)$ has a constant sectional curvature $\kappa$ such that
$R(X,Y,Z)=\kappa (g(Y,Z)X-g(X,Z)Y),\,\,X,Y,Z\in M_p$.
\end{itemize}
\end{Corollary}

$\bf{(b)\Longrightarrow (c)}$ Let $\{e_1, e_2, e_3, e_4\}$ be an
orthonormal basis for $M_p$, $p\in M$. Then the curvature operator
matrices $\mathcal{J}_{\{e_1,e_2\}}$ and
$\mathcal{J}_{\{e_3,e_4\}}$ have the form:

\begin{equation}
\mathcal{J}_{\{e_1,e_2\}}=
  \begin{pmatrix}
  K_{12}+K_{13} & 0 & R_{1223} & R_{1224}\\
  0 & K_{12} & R_{2113} & R_{2114}\\
  R_{1223} & R_{2113} & K_{13}+ K_{23} & \rho_{34}\\
  R_{1224} & R_{2114} & \rho_{34} & K_{14}+K_{24}
  \end{pmatrix},
 \end{equation}
and
\begin{equation}
\mathcal{J}_{\{e_3,e_4\}}=
  \begin{pmatrix}
    K_{13}+K_{14} & \rho_{12} & R_{1443} & R_{1334}\\
    \rho_{12} & K_{23}+K_{24} & R_{2443} & R_{2334}\\
    R_{1443} & R_{2443}& K_{34} & 0\\
    R_{1334} & R_{2334} & 0 & K_{34}
    \end{pmatrix}.
  \end{equation}

Let also

$$
(A)=\begin{pmatrix}
  K_{12}+K_{13} & 0 & R_{1223} & R_{1224}\\
  0 & K_{12} & R_{2113} & R_{2114}\\
  R_{1223} & R_{2113} & K_{13}+ K_{23} & \rho_{34}\\
  R_{1224} & R_{2114} & \rho_{34} & K_{14}+K_{24}
  \end{pmatrix}
\begin{pmatrix}
    K_{13}+K_{14} & \rho_{12} & R_{1443} & R_{1334}\\
    \rho_{12} & K_{23}+K_{24} & R_{2443} & R_{2334}\\
    R_{1443} & R_{2443}& K_{34} & 0\\
    R_{1334} & R_{2334} & 0 & K_{34}
    \end{pmatrix}.
$$
\indent Simple computations for the matrix $(A)$ give us:
\begin{align}
\hspace*{2cm}a_{11}&=K_{12}(K_{13}+K_{14})+R_{1223}R_{1443}+R_{1334}R_{1224},\nonumber\\
a_{12}&=K_{12}\rho_{12}+R_{1223}R_{2443}+R_{2334}R_{1224},\nonumber\\
a_{13}&=K_{12}R_{1443}+K_{34}R_{1223},\nonumber\\
a_{14}&=K_{12}R_{1334}+K_{34}R_{1224},\nonumber\\
a_{21}&=K_{12}\rho_{12}+R_{1443}R_{2113}+R_{2114}R_{1334},\nonumber \\
a_{22}&=K_{12}(K_{23}+K_{24})+R_{2443}R_{2113}+R_{2114}R_{2334},\nonumber\\
a_{23}&=K_{12}R_{2443}+K_{34}R_{2114},\nonumber\\
a_{24}&=K_{12}R_{2334}+K_{34}R_{2114},\nonumber\\
a_{31}&=(K_{13}+K_{13})R_{1223}+(K_{13}+K_{23})R_{1443}+\rho_{12}R_{2113}+\rho_{34}R_{1443},\\
a_{32}&=(K_{23}+K_{24})R_{2113}+(K_{13}+K_{23})R_{2443}+\rho_{12}R_{1223}+\rho_{34}R_{2334},\nonumber\\
a_{33}&=(K_{13}+K_{23})K_{34}+R_{1223}R_{1443}+R_{2443}R_{2113},\nonumber\\
a_{34}&=\rho_{34}K_{34}+R_{1223}R_{1334}+R_{2334}R_{2113},\nonumber\\
a_{41}&=(K_{13}+K_{14})R_{1224}+(K_{14}+K_{24})R_{1334}+\rho_{12}R_{2114}+\rho_{34}R_{1443},\nonumber\\
a_{42}&=(K_{23}+K_{24})R_{2114}+(R_{14}+K_{24})R_{2334}+\rho_{12}R_{1224}+\rho_{34}R_{2443},\nonumber\\
a_{43}&=\rho_{34}K_{34}+R_{1443}R_{1224}+R_{2114}R_{2334},\nonumber\\
a_{44}&=(K_{14}+K_{24})K_{34}+R_{1334}R_{1224}+R_{2114}R_{2334}.\nonumber
\end{align}
On the other hand let

$$
(B)=\begin{pmatrix}
  K_{13}+K_{14} & \rho_{12} & R_{1443} & R_{1334}\\
  \rho_{12} & K_{23}+K_{24} & R_{2443} & R_{2334}\\
  R_{1443} & R_{2443} & K_{34} & 0\\
  R_{1334} & R_{2334} & 0 & K_{34}
  \end{pmatrix}
\begin{pmatrix}
    K_{12} & 0 & R_{1223} & R_{1224}\\
    0 & K_{12} & R_{2113} & R_{2114}\\
    R_{1223} & R_{2113}& K_{13}+K_{23} & \rho_{34}\\
    R_{1224} & R_{2114} & \rho_{34} & K_{14}+K_{24}
    \end{pmatrix}.
$$

For the matrix $(B)$ elements we derive:

\begin{align}
\hspace*{2cm}b_{11}&=K_{12}(K_{13}+K_{14})+R_{1223}R_{1443}+R_{1334}R_{1224},\nonumber\\
b_{12}&=K_{12}\rho_{12}+R_{1443}R_{2113}+R_{2114}R_{1334},\nonumber\\
b_{13}&=(K_{13}+K_{14})R_{1223}+(K_{13}+K_{23})R_{1443}+\rho_{12}R_{2113}+\rho_{34}R_{1334},\nonumber\\
b_{14}&=(K_{13}+K_{14})R_{1224}+(K_{14}+K_{24})R_{1334}+\rho_{12}R_{2114}+\rho_{34}R_{1443},\nonumber\\
b_{21}&=K_{12}\rho_{12}+R_{1223}R_{2443}+R_{2334}R_{1224},\nonumber \\
b_{22}&=K_{12}(K_{23}+K_{24})+R_{2443}R_{2113}+R_{2114}R_{2334},\nonumber\\
b_{23}&=(K_{23}+K_{24})R_{2113}+(K_{13}+K_{23})R_{2443}+\rho_{12}R_{1223}+\rho_{34}R_{2334},\nonumber\\
b_{24}&= (K_{23}+K_{24})R_{2114}+(K_{14}+K_{24})R_{2334}+\rho_{12}R_{1224}+\rho_{34}R_{2443},\nonumber\\
b_{31}&=K_{12}R_{1443}+K_{34}R_{1223},\\
b_{32}&= K_{12}R_{2443}+K_{34}R_{2113},\nonumber\\
b_{33}&=(K_{13}+K_{23})K_{34}+R_{1223}R_{1443}+R_{2443}R_{2113},\nonumber\\
b_{34}&=\rho_{34}K_{34}+R_{1443}R_{1224}+R_{2114}R_{2443},\nonumber\\
b_{41}&= K_{12}R_{1334}+K_{34}R_{1224},\nonumber\\
b_{42}&= K_{12}R_{2334}+K_{34}R_{2114},\nonumber\\
b_{43}&=\rho_{34}K_{34}+R_{1223}R_{1334}+R_{2334}R_{2113},\nonumber\\
b_{44}&=(K_{14}+K_{24})K_{34}+R_{1334}R_{1224}+R_{2114}R_{2334}.\nonumber
\end{align}

From {\bf{[C2]}} we have $a_{12}=b_{12}$ and according to (2.6)
and (2.7), we get:
\medskip
\begin{itemize}
  \item[($e_7)$]$\quad$
  $R_{1223}R_{2443}+R_{2334}R_{1224}-R_{1443}R_{2113}-R_{2114}R_{1334}=0$,
  \item[$(e_8)$]$\quad$
 $ R_{1223}R_{1334}+R_{2334}R_{2113}-R_{1443}R_{1224}-R_{2114}R_{2443}=0$.
  \end{itemize}
\medskip
\indent By a cycling change of indecies $1\To 2\To 3\To 4\To
  1$ in $(e_7)$ and $(e_8)$ , we get
\medskip
\begin{itemize}
  \item[($e_7^1)$]$\quad$
  $R_{2334}R_{3114}+R_{1443}R_{1332}-R_{2114}R_{3224}-R_{1223}R_{1442}=0$,
  \item[$(e_8^1)$]$\quad$
 $ R_{2334}R_{1442}+R_{1443}R_{3224}-R_{2114}R_{1332}-R_{1223}R_{3114}=0$.
  \end{itemize}

We solve, using \textit{Maple}, $(e_7), (e_8), (e_7^1)$ and
$(e_8^1)$ together and arrive at the homogeneous system
\begin{equation}
  \begin{pmatrix}
    R_{2443} & R_{1224} & -R_{2113} & -R_{1334}\\
    R_{1334} & R_{2113} & -R_{1224} & -R_{2443}\\
    -R_{1442} & R_{3114} & R_{1332} & -R_{3224}\\
    -R_{3114} & R_{1442} & R_{3224} & -R_{1332}
    \end{pmatrix}
    \begin{pmatrix}
    R_{1223}\\
    R_{2334}\\
    R_{1443}\\
    R_{2114}
    \end{pmatrix}
    =
    \begin{pmatrix}
    0\\
    0\\
    0\\
    0
    \end{pmatrix}.
  \end{equation}

  Solving (2.8) with respect to $R_{1223}, R_{2334}, R_{1443}$ and
  $R_{2114}$, we get at least the trivial solution:

  \begin{equation}
  R_{2114}=R_{1223}=R_{1443}=R_{2334}=0.
  \end{equation}

  \indent We also have
  $\mathcal{J}_{\{e_1,e_3\}}=\mathcal{J}_{\{e_2,e_4\}}$ and
  $\mathcal{J}_{\{e_1,e_4\}}=\mathcal{J}_{\{e_2,e_3\}}$, and using
  (2.9), it follows that

  \begin{equation}
  \left|
\begin{array}{l}
R_{2113}(R_{1332}+R_{1442})+R_{1334}(R_{3114}+R_{3224})=0\\
 R_{1334}(K_{12}-K_{14}-K_{24})+R_{1224}(K_{34}-K_{13}-K_{14})=0\\
R_{2443}(K_{12}-K_{13}-K_{23})+R_{2113}(K_{34}-K_{23}-K_{24})=0\\
R_{1224}(R_{1332}+R_{1442})+R_{2443}(R_{3114}+R_{3224})=0\\
R_{3224}(R_{2113}+R_{2443})+R_{1442}(R_{1224}+R_{1334})=0\\
R_{1442}(K_{23}-K_{12}-K_{13})+R_{1332}(K_{14}-K_{34}-K_{13})=0\\
R_{1332}(R_{2113}+R_{2443})+R_{3114}(R_{1224}+R_{1334})=0
 \end{array}. \right.
 \end{equation}

 We solve (2.10) with respect to the tensor $R$ components
 $R_{1332}, R_{1442}, R_{3114}, R_{3224}$, $ R_{1224}, R_{1334},
 R_{2113}, R_{2114}$ using, for example \textit{Maple}, and as a
 result we get
 \begin{equation}
   \begin{array}{c}
    R_{1224}(-K_{34}+K_{13}+K_{14})=R_{1334}(K_{12}-K_{14}-K_{24})\\
    R_{2113}(-K_{34}+K_{23}+K_{24})=R_{2443}(K_{12}-K_{13}-K_{23})
    \end{array}
 \end{equation}

 and

 \begin{equation}
 R_{1332}=R_{3224}=R_{1442}=R_{3114}=0.
 \end{equation}
 \indent Further, by changing the orthonormal basis $\{e_1,e_2,
 e_3, e_4\}$ with the orthonormal basis $\left
 \{\dfrac{e_1-e_2}{\sqrt{2}},\dfrac{e_1+e_2}{\sqrt{2}},\dfrac{e_3-e_4}{\sqrt{2}},\dfrac{e_3+e_4}{\sqrt{2}}\right
 \}$ and using (2.11) and (2.12) we receive a new system of
 equations with respect to the tensor curvature components which
 is equivalent to (2.10). From there we can conclude that
 \begin{equation}
 R_{1334}=R_{2443}=0.
 \end{equation}
 \indent From (2.9), (2.11) and (2.12) it follows that all of the
 components $R_{ijjk}$ are equal to zero for all
 $i,j,k\,=1,2,3,4$.\\
 \indent If we reformulate the second and third equation in
 (2.10) by changing the basis as shown above and transforming
 the curvature components using (2.9), (2.11) and (2.12), we get
\begin{equation}
   \begin{array}{c}
    (K_{13}-K_{23}+K_{14}-K_{24})(R_{1432}+R_{1342}+K_{12})=0\\
    (K_{13}+K_{23}-K_{14}-K_{24})(R_{1432}+R_{1342}+K_{34})=0
    \end{array}.
 \end{equation}
 \indent By a cycling change of indecies $1\To 2\To 3\To 4\To 1$
 above and some extra tedious computations, we get the system:
  \begin{equation}
  \left|
\begin{array}{c}
(K_{13}+K_{14}-K_{23}-K_{24})(K_{12}+R_{1342}-R_{1423})=0\\
(K_{12}+K_{14}-K_{23}-K_{34})(K_{13}+R_{1234}-R_{1423})=0\\
(K_{12}+K_{13}-K_{24}-K_{34})(K_{14}+R_{1234}-R_{1342})=0\\
(K_{12}+K_{24}-K_{13}-K_{34})(K_{23}+R_{1234}-R_{1342})=0\\
(K_{12}+K_{23}-K_{14}-K_{34})(K_{24}+R_{1234}-R_{1423})=0\\
(K_{13}+K_{23}-K_{14}-K_{24})(K_{34}+R_{1342}-R_{1432})=0
 \end{array}. \right.
 \end{equation}

 Solving (2.14) and (2.15) with respect to the sectional curvature
 components $K_{12}, K_{13}$ and $K_{14}$, it follows that
 \begin{equation}
K_{14}=K_{23},\quad K_{13}=K_{24},\quad K_{12}=K_{34},
\end{equation}
 and since
 the basis $\{e_1, e_2, e_3, e_4\}$ has been arbitrary chosen in
 $M_p$, it follows that $(M,g)$ is Einstein \cite{Sin}.
 $\hfill\square$\\\medskip

 $\bf{(b)\Longrightarrow (a)}$ If $(M,g)$ is an Einstein manifold
  then $\rho=\lambda\Id$, $\lambda=\mbox{const.}$, and if
  $\alpha$ is a $2$-plane in $M_p,\,p\in M$, it follows that
  $$
 \begin{array}{c}
 \mathcal{J}(\alpha )\circ\mathcal{J}(\alpha^{\perp})-\mathcal{J}(\alpha^{\perp})\circ\mathcal{J}(\alpha )=
\mathcal{J}(\alpha
)\circ\mathcal{J}(\alpha^{\perp})+\mathcal{J}(\alpha
)\circ\mathcal{J}(\alpha )-\mathcal{J}(\alpha
)\circ\mathcal{J}(\alpha )
-\mathcal{J}(\alpha^{\perp})\circ\mathcal{J}(\alpha )=\\
\\
\mathcal{J}(\alpha )\circ\underbrace{\left
[\mathcal{J}(\alpha^{\perp})+\mathcal{J}(\alpha )\right ]}_{\rho}-
\underbrace{\left [\mathcal{J}(\alpha
)+\mathcal{J}(\alpha^{\perp})\right
]}_{\rho}\circ\mathcal{J}(\alpha )=\mathcal{J}(\alpha )\circ\rho
-\rho\circ\mathcal{J}(\alpha )=\\
\lambda\left (\mathcal{J}(\alpha
)\circ\Id-\Id\circ\mathcal{J}(\alpha )\right
)=0.\\
\end{array}
$$
   That completes the proof. $\hfill\square$

 \section{New approaches and results}

 Recently Gilkey, Puffini and Videv \cite{Vid4} were able to
 generalize the results above. They define $\mathfrak{M}:=(V,\langle\cdot,\cdot\rangle,A)$
 to be a $0$-model if
$\langle\cdot,\cdot\rangle$ is a non-degenerate inner product of
signature $(p,q)$ on a finite dimensional vector space $V$ of
dimension $m=p+q$ and if $A\in\otimes^4V^*$ is an algebraic
curvature tensor.\\
\indent Let $\Gr_{r,s}(V,\langle\cdot,\cdot\rangle)$ be the
Grassmannian of all non-degenerate linear subspaces of $V$ which
have  signature $(r,s)$; the pair $(r,s)$ is said to be {\it
admissible} if and only if $Gr_{r,s}(V,\langle\cdot,\cdot\rangle)$
is non-empty and does not consist of a single point or,
equivalently, if the inequalities $0\le r\le p$, $0\le s\le q$,
and $1\le r+s\le m-1$ are satisfied. Let $[A,B]:=AB-BA$ denote the
commutator of two linear maps. Then they establish the following
result:
\begin{thm}
Let $\mathfrak{M}=(V,\langle\cdot,\cdot\rangle,A)$ be a $0$-model.
The following assertions are equivalent; if any is satisfied, then
we shall say that $\mathfrak{M}$ is a {\rm Puffini--Videv}
$0$-model.
\begin{enumerate}
\item There exists $(r_0,s_0)$ admissible so that
\newline\qquad$\mathcal{J}(\pi)\circ\mathcal{J}(\pi^\perp)=
 \mathcal{J}(\pi^\perp)\circ\mathcal{J}(\pi)$ for all $\pi\in Gr_{r_0,s_0}(V,\langle\cdot,\cdot\rangle)$.
\item $\mathcal{J}(\pi)\circ\mathcal{J}(\pi^\perp)=
 \mathcal{J}(\pi^\perp)\circ\mathcal{J}(\pi)$ for every non-degenerate subspace $\pi$.
\item $[\mathcal{J}(\pi),\rho]=0$ for every non-degenerate
subspace $\pi$.
\end{enumerate}
\end{thm}

We say that $\mathfrak{M}=(V,\langle\cdot,\cdot\rangle,A)$ is {\it
decomposible} if there exists a non-trivial orthogonal
decomposition $V=V_1\oplus V_2$ which decomposes $A=A_1\oplus
A_2$; in this setting, we shall write
$\mathfrak{M}=\mathfrak{M}_1\oplus\mathfrak{M}_2$ where
$\mathfrak{M}_i:=(V_i,\langle\cdot,\cdot\rangle|_{V_i},A_i)$. One
says that $\mathfrak{M}$ is {\it indecomposible} if $\mathfrak{M}$
is not decomposible.

By Theorem 3.1, any Einstein $0$-model is Puffini--Videv. More
generally, the direct sum of Einstein Puffini--Videv models is
again Puffini--Videv; the converse holds in the Riemannian
setting:
\begin{thm}
Let $\mathfrak{M}=(V,\langle\cdot,\cdot\rangle,A)$ be a Riemannian
$0$-model. Then $\mathfrak{M}$ is Puffini--Videv if and only if
$\mathfrak{M}=\mathfrak{M}_1\oplus\cdots\oplus\mathfrak{M}_k$
where the $\mathfrak{M}_i$ are Einstein.
\end{thm}

In the pseudo-Riemannian setting, a somewhat weaker result can be
established. One says that a $0$-model is {\it pseudo-Einstein}
either if the Ricci operator $\rho$ has only one real eigenvalue
$\lambda$ or if the Ricci operator $\rho$ has two complex
eigenvalues $\lambda_1,\lambda_2$ with $\bar\lambda_1=\lambda_2$.
This does not imply that $\rho$ is diagonalizable in the higher
signature setting and hence $\mathfrak{M}$ need not be Einstein.

\begin{thm}
 Let $\mathfrak{M}=(V,\langle\cdot,\cdot\rangle,A)$ be a $0$-model of arbitrary signature. If $\mathfrak{M}$ is Puffini--Videv,
 then we may decompose
$\mathfrak{M}=\mathfrak{M}_1\oplus\cdots\oplus\mathfrak{M}_k$ as
the direct sum of pseudo-Einstein $0$-models $\mathfrak{M}_i$.
\end{thm}

\bibliographystyle{plain}
\nocite*
\bibliography{Article_PU}

\begin{thebibliography}{10}

\bibitem{Bro}
M.~Brozos-V{\'a}zquez and P.~Gilkey.
\newblock The global geometry of {R}iemannian manifolds with comuting curvature
  operators.
\newblock preprint.

\bibitem{Chi}
Q.-Sh. Chi.
\newblock A curvature characterization of certain locally rank-one symmetric
  spaces.
\newblock {\em J. Diff. Geom.}, 28:187--202, 1988.

\bibitem{Rio}
E.~Garc{\'i}a-R{\'i}o, D.~Kupeli, and R.~V{\'a}zquez-Lorenzo.
\newblock {\em Osserman {M}anifolds in {S}emi-{R}iemannian {G}eometry}.
\newblock Number 1777 in Lecture Notes in Math. Springer-Verlag, Berlin, 2002.

\bibitem{Vid4}
P.~Gilkey, E.~Puffini, and V.~Videv.
\newblock Puffini-{V}idev {M}odels and {M}anifolds.
\newblock {\em J. of Geom. (to appear)}, 2006.
\newblock ar{X}iv: math.DG/0605464.

\bibitem{Vid1}
P.~Gilkey, G.~Stanilov, and V.~Videv.
\newblock Pseudo-{R}iemannian manifolds whose generalized {J}acobi operator has
  a constant characteristic polynomial.
\newblock {\em J. of Geom.}, 62:144--153, 1998.

\bibitem{Kob}
S.~Kobayashi and K.~Nomizu.
\newblock {\em Foundations of Differential Geometry}, volume~I.
\newblock Interscience Publish., New York, 1963.

\bibitem{Nik}
Y.~Nikolaevsky.
\newblock Two theorems on {O}sserman manifolds.
\newblock {\em Diff. Geom. {\&} Appl.}, 18:239--253, 2003.

\bibitem{Oss}
R.~Osserman.
\newblock Curvature in the eighties.
\newblock {\em Amer. Math. Montly}, 97:731--756, 1990.

\bibitem{Schouten}
J.A. Schouten and D.~J. Struik.
\newblock On some properties of general manifolds relating to {E}instein's
  theory of gravitation.
\newblock {\em Amer. J. Math}, 43:213--216, 1921.

\bibitem{Sin}
I.~Singer and J.~Thorpe.
\newblock The curvature of 4-dimensional {E}instein spaces.
\newblock {\em Global Analysis}, 1969.
\newblock Papers in Honor of K. Kodaira. University of Tokyo Press, Princeton
  University Press.

\bibitem{Vid2}
G.~Stanilov and V.~Videv.
\newblock On the commuting {S}tanilov's curvature operators.
\newblock {\em Mathematics and Education in Mathematics (Proc. of the 33rd
  {S}pring {C}onference of the {U}nion of {B}ulgarian {M}athematicians,
  Borovets, April 1-4, 2004)}, pages 176--180, 2004.

\bibitem{Stan1}
Gr. Stanilov.
\newblock {\em Differential geometry}.
\newblock Nauka i izkustvo, Sofia, 1988.
\newblock (in Bulgarian).

\bibitem{Vid3}
V.~Videv and M.~Ivanova.
\newblock Four-dimensional {R}iemannian manifolds with commuting {S}tanilov
  curvature operators with respect to the orthogonal planes.
\newblock {\em Mathematics and Education in Mathematics (Proc. of the 33rd
  {S}pring {C}onference of the {U}nion of {B}ulgarian {M}athematicians,
  Borovets, April 1-4, 2004)}, pages 180--184, 2004.

\bibitem{Wang}
H.~C. Wang.
\newblock Two-point homogeneous spaces.
\newblock {\em Ann. of Math.}, 55:177--191, 1952.

\end{thebibliography}

\end{document}